%Paper: alg-geom/9411011
%From: miranda@riemann.MATH.ColoState.EDU (Rick Miranda)
%Date: Wed, 16 Nov 94 14:58:12 MST

\magnification=1200
\baselineskip 18pt
\def\I{I\!\!P}
\font\ms=msbm10

\font\msq=msbm7
\def\C{\ms C}

\def\qed{\vrule height 1.1ex width 1.0ex depth -.1ex }
\def\mapright#1{\smash{\mathop{\longrightarrow}\limits^{#1}}}

\centerline{\bf CLASSIFICATION OF VARIETIES WITH CANONICAL CURVE
SECTION}
\centerline{\bf VIA GAUSSIAN MAPS ON CANONICAL CURVES}
\vskip 1.5cm
\centerline{CIRO CILIBERTO*\footnote\null{{\it \noindent 1991
Mathematics Subject Classification:} Primary 14J45. Secondary
14J28, 14H10.},
ANGELO FELICE LOPEZ\footnote*{Research  partially supported by
the MURST
national project ``Geometria Algebrica";} AND
RICK\footnote\null{the authors are
members of GNSAGA of CNR.} MIRANDA*\footnote*{* Research
partially supported by
NSF under grant DMS-9403392.}}
\vskip 1.5cm
\noindent {\bf 1. INTRODUCTION}
\vskip .5cm

Let $C \subset \I^{g-1}$ be a smooth canonical curve of genus
$g \geq 3$.
The purpose of this article is to further develop a method to
classify varieties
having $C$ as their curve section, following the techniques of
[CLM1]. In that
article a careful analysis of the degeneration to the cone over
the hyperplane
section was made for {\it prime} Fano threefolds, that is Fano
threefolds whose
Picard group is generated by the hyperplane bundle. In this
article we extend
this method and classify Fano threefolds of higher index (which
still have
Picard number one) (see Theorem (3.2)). We are also able to
classify Mukai
varieties, i.e. varieties of dimension four or more with canonical
curve
sections (see Theorem (3.11)).

The classification of the varieties in question has been executed
before, mostly
by Fano, with modern proofs given by Iskovskih, Mukai and others.
These proofs
rely on several deep theorems, in particular for the Fano
threefolds, the
existence of lines and of smooth sections of in the primitive
submultiple of the
anticanonical line bundle.

Our approach to the classification is completely different. It
is based on a
general remark concerning the degeneration to the cone over the
hyperplane
section, as mentioned above. Let $V\subset \I^r$ be a smooth,
irreducible,
projectively Cohen-Macaulay variety, with general hyperplane
section $W$. Note
that $V$ flatly degenerates to the cone $X$ over $W$, and therefore
$[X]$ is a
point in the Hilbert scheme $\cal H$ of $V$. Suppose one has an
upper bound for
$h^0(X,N_X)$, which is the dimension of the Zariski tangent space
to $\cal H$ at
$[X]$. Suppose in addition that one has a known irreducible family
$\cal F$ of
varieties containing $[V]$ as a member, and that the dimension of
$\cal F$ is
equal to the upper bound for $h^0(X,N_X)$. Since the closure of
$\cal F$ in the
Hilbert scheme $\cal H$ must contain the point $[X]$, we conclude
that the
dimension of $\cal H$ at $[X]$ is equal to the dimension of its
Zariski tangent
space there, and hence $[X]$ is a smooth point of $\cal H$, and
the closure of
the known family $\cal F$ is the unique component of $\cal H$
containing $[X]$
(and $[V]$).

One more ingredient now comes into play. Suppose further that the
Hilbert scheme
$\cal H'$ of the hyperplane sections $W$ of $V$ is irreducible, and
that one can
prove that as $V$ varies in any component of its Hilbert scheme, the
possible
hyperplane sections $W$ vary to fill up $\cal H'$. This assumption
allows us to
conclude then that the Hilbert scheme $\cal H$ is also irreducible:
if there
were two components, each would contain varieties with general
hyperplane
section $W$, and therefore each would contain the cone point $[X]$
over a
general $W$. This would force $[X]$ to be a singular point of the
Hilbert
scheme, which, given the first part of the argument as described
above, is a
contradiction. One concludes then that there is only one component of
$\cal H$,
which is the closure of the known family $\cal F$ as described above.
In this
sense one obtains a classification result for varieties $V$.

This scheme can be applied to classify Fano threefolds of the principal
series
with Picard number one. First of all from papers of Fano and Iskovskih
one may
easily extract lists of families of such threefolds, and the number of
parameters on which they depend are easily available. (We stress that
our
classification argument does not rely on the classification theorems of
Fano and
Iskovskih, but simply on the existence of the families, taken as
examples.)
Secondly an easy deformation theoretic argument proves that the general
hyperplane section of a general Fano threefold of the principal series
varying
in any component of the Hilbert scheme, is a general K3 surface of the
correct
genus. Therefore to complete the program outlined above, we must give
the sharp
upper bound for $h^0(X,N_X)$. By now it is well known that this
computation rests
on the study of Gaussian maps for the general curve section $C$. This
is done
in section 2 (see Theorem (2.13)).

We are also able to execute the program in an iterative fashion, to
obtain
similar classification statements for the Mukai varieties. The list
of examples
is also available, together with the number of parameters on which
they depend.

The reader will notice that the above scheme provides a
classification
statement only when there are actually examples of such varieties.
To finish one
has to prove a negative statement, that varieties do not exist when
there are no
examples. This is provided by a theorem of Zak and L'vovskii, which,
under
suitable conditions on the Gaussian map of a curve, says that this
curve cannot
be the curve section a variety which is not a cone. The Gaussian map
computations which we have to perform to prove the classification part,
suffice
also to demonstrate the required non existence statement for the Fano
and Mukai
varieties.

In section 2 we will recall the definition of the relevant Gaussian maps
and we
will perform the computations of the coranks of these maps which are
necessary
for our purposes, as indicated above. In section 3 we execute the
classification
program for the Fano threefolds and Mukai varieties. Since the case
of prime
Fano threefolds was already treated in [CLM1], we treat here only the
cases of
higher index.
\vskip .5cm

\noindent {\bf 2. THE CORANK THEOREM}

\vskip .5cm

We recall here
the definition of the {\it Gaussian map}
$$\matrix{\Phi_{\omega_{C}} : \bigwedge^2 H^0(C, \omega_{C}) \to
H^0(C,\omega_{C}^{\otimes 3}) \cr}$$
which sends $s\wedge t$ to $s\otimes dt-t \otimes ds$. The Corank
One Theorem
from [CLM1]
(Theorem 4) states that for a general hyperplane section $C_g$ of a
general
{\it prime} K3 surface $S_g$ (that is one whose Picard group is
generated by the
hyperplane bundle) the corank of the Gaussian map $\Phi_{\omega_{C_g}}$
is one
if $g = 11$ or $g \geq 13$.

In this section we will extend this result to {\it reembedded} K3
surfaces.
Let ${\cal H}_{g}$ be the component of the Hilbert scheme of prime K3
surfaces of
genus $g \geq 3$. Let ${\cal H}_{2}$ be the family of prime genus $2$
K3
surfaces, that is double covers of $\I^2$ ramified along a sextic curve.
For $r
\geq 2$ let ${\cal H}_{r,g}$ be the component of the Hilbert scheme
whose
general elements are obtained by embedding prime K3 surfaces of genus
$g$ via
the $r$-th multiple of the primitive class, except for $r = g = 2$ when
we set
${\cal H}_{2,2} = {\cal H}_{2}$. Moreover set ${\cal H}_{1,g} =
{\cal H}_{g}$.
We denote by $S_{r,g}$ any smooth surface represented by a point
in ${\cal
H}_{r,g}$; note that, for $(r,g) \not= (2,2)$, $S_{r,g} \subset
\I^{N(r,g)}$,
where $N(r,g) = 1 + r^2(g-1)$. We will say that $S_{r,g}$ is {\it non
trigonal} if there is no trigonal curve in its primitive linear system.
We
finally denote by $C_{r,g}$ any smooth curve section of $S_{r,g}$, if
$(r,g)
\not= (2,2)$, and by $C_{2,2}$ any smooth curve in the linear system on
$S_{2,2}$ given by twice the primitive linear system.

The following general remark is the basis of our computation. Let $C$ be
a smooth hyperplane
section of a smooth K3 surface $S$ and consider the commutative diagram
$$\matrix{\bigwedge^2 H^0(S,{\cal O}_{S}(1))  &
\mapright{\Phi_{{\cal O}_{S}(1)}} & H^0(S,\Omega_{S}^1(2)) \cr
&& \downarrow
\phi \cr  \downarrow \pi & & \ \ H^0(C,\Omega_{S}^1(2)_{|C}) \cr
&&\downarrow
\psi \cr \bigwedge^2 H^0(C,{\cal O}_{C}(1))
&\mapright{\Phi_{\omega_{C}}}  &
H^0(C,\omega_{C}^{\otimes 3}) \cr} \leqno (2.1)$$

\noindent where $\pi$ and $\phi$ are restriction maps,
$\Phi_{{\cal O}_{S}(1)}$ and $\Phi_{\omega_{C}}$ are Gaussian maps and
$\psi$ is the map on differentials induced by the sheaf map
${\Omega_{S}^1}_{|C} \to \Omega_{C}^1$. We have the following:

{\bf Lemma (2.2).} {\sl Let $C$ be a smooth hyperplane section of a
smooth K3
surface $S$, $\phi$ and $\psi$ as in diagram (2.1). If
$\Phi_{{\cal O}_{S}(1)}$
and $\phi$ are surjective then $\Phi_{\omega_{C}}$ has corank one.}

\noindent {\it Proof:} From (2.1) we have $Im \ \psi = Im \ (\psi
\circ \phi
\circ \Phi_{{\cal O}_{S}(1)}) = Im \ (\Phi_{\omega_{C}} \circ \pi)
= Im
\ \Phi_{\omega_{C}}$ since $\phi \circ \Phi_{{\cal O}_{S}(1)}$ and
$\pi$ are
surjective. On the other hand the exact sequence
$$ 0 \to N_{C/S}^*(2) \to \Omega_{S}^1(2)_{|C} \to
\omega_{C}^{\otimes 3} \to 0$$

\noindent shows that $corank \ \psi \leq h^1(N_{C/S}^*(2)) =
h^1({\cal
O}_{C}(1)) = h^1(\omega_{C}) = 1,$ hence $corank \
\Phi_{\omega_{C}} = 1$ since
$\Phi_{\omega_{C}}$ is not surjective by Wahl's theorem ([W])
because $C$ lies
on a K3 surface. \ \qed

By (2.1) and Lemma (2.2) the computation of the
corank of the Gaussian map  $\Phi_{\omega_{C_{r,g}}}$ is then
reduced to
studying the maps $\Phi_{{\cal O}_{S_{r,g}}(1)}$ and  $\phi_{r,g}:
H^0(S_{r,g},\Omega_{S_{r,g}}^1(2)) \to
H^0(S_{r,g},\Omega_{S_{r,g}}^1(2)_{|C_{r,g}})$. The corank of
$\Phi_{{\cal O}_{S_{r,g}}(1)}$ was calculated in [CLM2,
Theorem (1.1)]. As for
$\phi_{r,g}$ we have:

{\bf Lemma (2.3).} {\sl Suppose $k \geq 1$ and $r, g$ are such that
either
$r = 1, g = 11$ or $g \geq 13$; or $ r \geq 2, g \geq 2$. Then
$H^1(S_{r,g},\Omega_{S_{r,g}}^1(k)) = 0$ and therefore $\phi_{r,g,k}:
H^0(S_{r,g},\Omega_{S_{r,g}}^1(k+1)) \to
H^0(S_{r,g},\Omega_{S_{r,g}}^1(k+1)_{|C_{r,g}})$ is surjective if one
of the
following holds:

(a) $k \geq 2$ and $r \geq 2, g \geq 3$ except $k = r = 2, g = 3$;

(b) $k = 1$ and $r = 4, g \geq 4$ or $r \geq 5, g \geq 3$;

(c) $S_{r,g}$ represents a general point of ${\cal H}_{r,g}$, $k = 1$
and
$r = 1, g = 11, g \geq 13$ or
\indent $r = 2, g \geq 7$;

(d) $S_{3,g}$ is non trigonal, $k = 1, r = 3$ and $g \geq 5$;

(e) $S_{r,2}$ has smooth ramification divisor, $g = 2, k = 1$ and $r
\geq 4$ or
$k \geq 2$ and $r \geq 3$.

\noindent Moreover for $k = r = 2, g = 3$ or $k = 1$ and $r = 4, g =
3$ or $r =
2, g = 6$ and  $S_{2,6}$ non trigonal, we have instead $corank \
\phi_{r,g,k} =
1$.}

\noindent {\it Proof:} By definition of $S_{r,g}$ we have
$$H^1(S_{r,g},\Omega_{S_{r,g}}^1(k)) \cong
H^1(S_{g},\Omega_{S_{g}}^1(rk)) \cong H^1(S_{g},T_{S_{g}}(-rk))^*$$
\noindent where $S_{g}$ is a smooth K3 surface representing a point
of
${\cal H}_{g}$.

\noindent Consider now for $g \geq 3$ the Euler sequence of
$S_{g} \subset \I^g$ tensored by ${\cal O}_{S_{g}}(-j)$
$$0 \to {\cal O}_{S_{g}}(-j) \to H^0({\cal O}_{S_{g}}(1))^*
\otimes  {\cal O}_{S_{g}}(1-j) \to {T_{\I^g}}_{|S_{g}}(-j) \to 0.
\leqno
(2.4)$$
\noindent Let us show that
$$h^0({T_{\I^g}}_{|S_{g}}(-j)) = \cases{ g+1 &for $j=1$ \cr $0$ &for
$j \geq
2$\cr} \leqno (2.5)$$
\noindent and
$$H^1({T_{\I^g}}_{|S_{g}}(-j)) = 0 \ \ \hbox{for} \ j \geq 1. \leqno
(2.6)$$
\noindent For $j \geq 1$, $H^0({\cal O}_{S_{g}}(-j)) = 0$ and
$H^1({\cal
O}_{S_{g}}(-j)) = 0$ by Kodaira vanishing, hence (2.4) gives
$h^0({T_{\I^g}}_{|S_{g}}(-j)) = (g+1)h^0({\cal O}_{S_{g}}(1-j))$,
which is
(2.5). Also $H^1({\cal O}_{S_{g}}(1-j)) = 0$ by Kodaira vanishing
for $j \geq 2$
and for $j = 1$ because $S_{g}$ is a K3 surface, therefore by (2.4)
we have
$$H^1({T_{\I^g}}_{|S_{g}}(-j)) = Ker \ \{ H^2({\cal O}_{S_{g}}(-j))
\to H^0({\cal
O}_{S_{g}}(1))^* \otimes  H^2({\cal O}_{S_{g}}(1-j)) \} = $$
$$ = ( Coker \ \{ H^0({\cal O}_{S_{g}}(1)) \otimes
H^0({\cal O}_{S_{g}}(j-1))
\to H^0({\cal O}_{S_{g}}(j)) \})^* = 0$$
\noindent by the surjectivity of this multiplication map, hence
(2.6). Now the
normal bundle sequence
$$0 \to T_{S_{g}}(-j) \to {T_{\I^g}}_{|S_{g}}(-j) \to
N_{S_{g}}(-j) \to 0$$
\noindent gives
$$h^1(T_{S_{g}}(-j)) = \cases{ h^0(N_{S_{g}}(-1)) - g -1 &for
$j=1$ \cr
h^0(N_{S_{g}}(-j)) &for $j \geq 2$\cr}. \leqno (2.7)$$
\noindent In fact for $j \geq  2$ we have
$H^0({T_{\I^g}}_{|S_{g}}(-j)) =
H^1({T_{\I^g}}_{|S_{g}}(-j)) = 0$ by (2.5) and (2.6),
while for $j = 1$, since $h^2(\Omega_{S_{g}}^1) =
h^{1,2}(S_{g}) = 0$ it is a fortiori $h^0(T_{S_{g}}(-1)) =
h^2(\Omega_{S_{g}}^1(1)) = 0$, hence (2.7) holds again by (2.5) and
(2.6).

Now let $C_{g}$ be a general hyperplane section of $S_{g}$, $X_{g}$
a
general cone over $C_{g}$ in $\I^g$. Then $S_{g}$ flatly
degenerates to
$X_{g}$ ([CLM1], 5.2) and $h^0(N_{S_{g}}(-j))$ is upper
semicontinuous, hence
$$h^0(N_{S_{g}}(-j)) \leq h^0(N_{X_{g}}(-j)) =
\sum\limits_{h \geq 0}
h^0(N_{C_{g}/\I^{g-1}}(-j-h)). \leqno (2.8)$$
\noindent We now set $j = rk$. For $k = r = 1, g = 11$ or $g
\geq 13$ and
$S_{g}$ general we have that
$$h^0(N_{C_{g}/\I^{g-1}}(-1-h)) = \cases{ g + 1
&for $h=0$ \cr $0$ &for $h \geq 1$ \cr}$$
\noindent by [CLM1, Theorem 4 and Lemma 4], hence
$h^0(N_{S_{g}}(-1)) \leq g + 1$ and
$h^1(S_{1,g},\Omega_{S_{1,g}}^1(1)) =
h^1(S_{g},T_{S_{g}}(-1)) = 0$ by (2.7). For $k = 1, r = 2, g
\geq 7$ and
$S_{g}$ general we have that
$h^0(N_{C_{g}/\I^{g-1}}(-2-h)) = 0$ for $h \geq 0$ again by [CLM1,
Lemma 4].
Hence by (2.7) and (2.8) we deduce
$h^1(S_{2,g},\Omega_{S_{2,g}}^1(1)) =
h^1(S_{g},T_{S_{g}}(-2)) = h^0(N_{S_{g}}(-2)) = 0$. This proves
the
vanisihing under hypothesis (c). Similarly, if we assume (d), that
is for $k =
1, r = 3, g \geq 5$ and $S_{g}$ non trigonal, we have that
$h^0(N_{C_{g}/\I^{g-1}}(-3-h)) = 0$ for $h \geq 0$ because the ideal
of
$C_{g}$ is generated by quadrics, hence
$h^1(S_{3,g},\Omega_{S_{3,g}}^1(1)) =
h^1(S_{g},T_{S_{g}}(-3)) = h^0(N_{S_{g}}(-3)) = 0$.

In hypothesis (e) we have $g = 2, k = 1$ and $r \geq 4$ or $k \geq 2,
r \geq 3$
and $S_{2}$ is a double cover $\pi : S_{2} \to \I^2$ with smooth
ramification
divisor $B$. The Euler sequence    $$0 \to {\cal O}_{S_{2}}(-rk)
\to H^0({\cal
O}_{S_{2}}(1))^* \otimes  {\cal O}_{S_{2}}(1-rk) \to
\pi^*T_{\I^2}(-rk) \to 0$$
\noindent gives, as above by
Kodaira vanishing, that $H^0(\pi^*T_{\I^2}(-rk)) =
H^1(\pi^*T_{\I^2}(-rk)) = 0$;
the normal bundle sequence is
$$0 \to T_{S_{2}}(-rk) \to \pi^*T_{\I^2}(-rk) \to
N_{\pi}(-rk) \to 0$$
\noindent and $N_{\pi} \cong \pi^*{\cal O}_{B}(3)$ and therefore we
get that
$h^1(S_{r,2},\Omega_{S_{r,2}}^1(k)) = h^1(S_{2},T_{S_{2}}(-rk)) =
h^0(N_{\pi}(-rk)) = h^0(\pi^*{\cal O}_{B}(3-rk)) = 0$. This gives
the required
vanishings under hypothesis (e). Suppose now that we are under
hypotheses (a)
or (b). Then we have
$h^0(N_{C_{g}/\I^{g-1}}(-rk-h)) = 0$ for $h \geq 0$ since the
ideal of $C_{g}$
is generated by hypersurfaces of degree $\leq rk - 1$ and
therefore
$h^1(S_{r,g},\Omega_{S_{r,g}}^1(k)) = h^1(S_{g},T_{S_{g}}(-rk)) =
h^0(N_{S_{g}}(-rk)) = 0$ by (2.7) and (2.8).

For $g = 3$ and $k = r = 2$ or $k = 1, r = 4$ we have that
$H^1(S_{r,3},\Omega_{S_{r,3}}^1(k+1)) = 0$ by part (a) and from
(2.7) we get
$h^1(S_{r,3},\Omega_{S_{r,3}}^1(k)) = h^1(S_{3},T_{S_{3}}(-rk)) =
h^0(N_{S_{3}}(-4)) = h^0({\cal O}_{S_{3}}) = 1$, hence $corank \
\phi_{r,3,k} =
1$. For $k = 1, r = 2, g = 6$ and $S_{2,6}$ non trigonal, we will
prove that
$h^1(S_{2,6},\Omega_{S_{2,6}}^1(1)) =
h^1(S_{6},\Omega_{S_{6}}^1(2)) = 1$,
hence $corank \ \phi_{2,6,1} = 1$ because
$H^1(S_{2,6},\Omega_{S_{2,6}}^1(2)) =
0$ by part (a). Let $\hbox{\ms G} = \hbox{\ms G}(1,4)$ be the
Grassmannian
of lines in $\I^4$ and ${\cal O}_{\hbox{\ms G}}(1)$ the line bundle
giving the
Pl\"ucker embedding. Recall that the K3 surface $S_{6}$ is
scheme-theoretically
cut out by three sections of ${\cal O}_{\hbox{\ms G}}(1)$ and one
section of
${\cal O}_{\hbox{\ms G}}(2)$. From the normal bundle sequence
$$0 \to N_{S_{6} /
\hbox{\msq G}}^*(2) \to  \Omega_{\hbox{\ms G}}^1(2)_{|S_{6}}
\to \Omega_{S_{6}}^1(2) \to 0$$
\noindent and the fact that
$h^2(N_{S_{6} / \hbox{\msq G}}^*(2)) =
h^0({\cal O}_{S_{6}}(-1)^{\oplus 3} \oplus {\cal O}_{S_{6}}) = 1$,
we are
reduced to prove
$$H^p(\Omega_{\hbox{\ms G}}^1(2)_{|S_{6}}) = 0 \ \ \hbox{for} \
p =1, 2. \leqno (2.9)$$
\noindent From the exact sequence
$$0 \to \Omega_{\hbox{\ms G}}^1 \otimes
{\cal I}_{S_{6} /\hbox{\msq G}}(2) \to \Omega_{\hbox{\ms G}}^1 (2)
\to
\Omega_{\hbox{\ms G}}^1(2)_{|S_{6}} \to 0$$
\noindent we see that (2.9) follows once we prove
$$H^p(\Omega_{\hbox{\ms G}}^1 \otimes
{\cal I}_{S_{6} /\hbox{\msq G}}(2)) = 0 \ \ \hbox{for} \ p = 2, 3
\leqno (2.10)$$
\noindent since $H^p(\Omega_{\hbox{\ms G}}^1(2)) = 0$ for $p = 1,
2$ by Bott
vanishing. But also (2.10) follows easily by Bott vanishing since the
Koszul
resolution of the ideal sheaf ${\cal I}_{S_{6} /\hbox{\msq G}}$ gives
the
exact sequence
$$0 \to {\Omega_{\hbox{\ms G}}^1}(-3) \to {\Omega_{\hbox{\ms
G}}^1}(-2)^{\oplus 3} \oplus {\Omega_{\hbox{\ms G}}^1}(-1) \to
{\Omega_{\hbox{\ms G}}^1}^{\oplus 3} \oplus {\Omega_{\hbox{\ms
G}}^1}(-1)^{\oplus 3} \to$$
$$ \to {\Omega_{\hbox{\ms G}}^1}(1)^{\oplus 3}
\oplus {\Omega_{\hbox{\ms G}}^1} \to {\Omega_{\hbox{\ms G}}^1} \otimes
{\cal I}_{S_{6} / \hbox{\msq G}}(2) \to 0.  \ \ \ \ \qed$$
{\bf (2.11)} {\it Remark.} Part of the above Lemma, namely the fact
that $H^1(S_{1,g},\Omega_{S_{1,g}}^1(1)) = H^1(S_{g},T_{S_{g}}(-1)) =
0$ for $g
= 11$ or $g \geq 13$ also follows by work of Mori and Mukai ([MM],
[M1]).
\eject
Before coming to the statements and proofs of the theorems announced
in the
introduction, we will collect in one single table all the useful
information
about the K3 surfaces $S_{r,g} \subset \I^{N(r,g)}$, their hyperplane
sections
$C_{r,g}$ and the cones ${\widehat S}_{r,g} \subset \I^{N(r,g) + 1}$
over
$S_{r,g}$. This will be proved and used along the paper:

\centerline{Table (2.12)}

\centerline{\vbox{\offinterlineskip
\hrule
\halign{&\vrule#&
\strut\hfil#\tabskip=.001em\cr
height4pt&\omit&&\omit&&\omit&&\omit&&\omit&&\omit&&\omit&\cr
&$r$\hfil&&$g$\hfil&&\ $corank \
\Phi_{\omega_{C_{r,g}}}$&&\ $h^0(N_{C_{r,g}}(-2))$&&\ bound on
$h^0(N_{{\widehat
S}_{r,g}})$&&\ hypothesis on $S_{r,g}$&&\ hypothesis on $C_{r,g}$&\cr
height4pt&\omit&&\omit&&\omit&&\omit&&\omit&&\omit&&\omit&\cr
\noalign{\hrule}
height4pt&\omit&&\omit&&\omit&&\omit&&\omit&&\omit&&\omit&\cr
&1\hfil&&6\hfil&&10\hfil&&$\leq
1$\hfil&&85\hfil&&general\hfil&&general\hfil&\cr
height4pt&\omit&&\omit&&\omit&&\omit&&\omit&&\omit&&\omit&\cr
&1\hfil&&7\hfil&&9\hfil&&0\hfil&&98\hfil&&general\hfil&&general
\hfil&\cr
height4pt&\omit&&\omit&&\omit&&\omit&&\omit&&\omit&&\omit&\cr
&1\hfil&&8\hfil&&7\hfil&&0\hfil&&114\hfil&&general\hfil&&general
\hfil&\cr
height4pt&\omit&&\omit&&\omit&&\omit&&\omit&&\omit&&\omit&\cr
&1\hfil&&9\hfil&&5\hfil&&0\hfil&&132\hfil&&general\hfil&&general
\hfil&\cr
height4pt&\omit&&\omit&&\omit&&\omit&&\omit&&\omit&&\omit&\cr
&1\hfil&&10\hfil&&4\hfil&&0\hfil&&153\hfil&&general\hfil&&general
\hfil&\cr
height4pt&\omit&&\omit&&\omit&&\omit&&\omit&&\omit&&\omit&\cr
&1\hfil&&11\hfil&&1\hfil&&0\hfil&&\null\hfil&&general\hfil&&general
\hfil&\cr
height4pt&\omit&&\omit&&\omit&&\omit&&\omit&&\omit&&\omit&\cr
&1\hfil&&12\hfil&&2\hfil&&0\hfil&&201\hfil&&general\hfil&&general
\hfil&\cr
height4pt&\omit&&\omit&&\omit&&\omit&&\omit&&\omit&&\omit&\cr
&1\hfil&&\
$\geq 13$\hfil&&1\hfil&&0\hfil&&\null\hfil&&general\hfil&&general
\hfil&\cr
height4pt&\omit&&\omit&&\omit&&\omit&&\omit&&\omit&&\omit&\cr
&1\hfil&&\
$\geq 17$\hfil&&1\hfil&&0\hfil&&\null\hfil&&general\hfil&&any
smooth\hfil&\cr
height4pt&\omit&&\omit&&\omit&&\omit&&\omit&&\omit&&\omit&\cr
&2\hfil&&2\hfil&&13\hfil&&\null\hfil&&\null\hfil&&any
smooth\hfil&&any
smooth\hfil&\cr
height4pt&\omit&&\omit&&\omit&&\omit&&\omit&&\omit&&\omit&\cr
&2\hfil&&3\hfil&&10\hfil&&$\leq
1$\hfil&&139\hfil&&any smooth\hfil&&any
smooth\hfil&\cr
height4pt&\omit&&\omit&&\omit&&\omit&&\omit&&\omit&&\omit&\cr
&2\hfil&&4\hfil&&7\hfil&&0\hfil&&234\hfil&&any smooth\hfil&&any
smooth\hfil&\cr
height4pt&\omit&&\omit&&\omit&&\omit&&\omit&&\omit&&\omit&\cr
&2\hfil&&5\hfil&&4\hfil&&0\hfil&&363\hfil&&non trigonal\hfil&&any
smooth\hfil&\cr
height4pt&\omit&&\omit&&\omit&&\omit&&\omit&&\omit&&\omit&\cr
&2\hfil&&6\hfil&&2\hfil&&0\hfil&&525\hfil&&non trigonal\hfil&&any
smooth\hfil&\cr
height4pt&\omit&&\omit&&\omit&&\omit&&\omit&&\omit&&\omit&\cr
&2\hfil&&\ $\geq 7$\hfil&&1\hfil&&0\hfil&&\null\hfil&&general
\hfil&&any
smooth\hfil&\cr
height4pt&\omit&&\omit&&\omit&&\omit&&\omit&&\omit&&\omit&\cr
&3\hfil&& 2\hfil&&10\hfil&&1\hfil&&\null\hfil&&smooth ramif.
div.\hfil&&general\hfil&\cr
height4pt&\omit&&\omit&&\omit&&\omit&&\omit&&\omit&&\omit&\cr
&3\hfil&&2\hfil&&18\hfil&&1\hfil&&\null\hfil&&smooth ramif.
div.\hfil&&special\hfil&\cr
height4pt&\omit&&\omit&&\omit&&\omit&&\omit&&\omit&&\omit&\cr
&3\hfil&&3\hfil&&5\hfil&&0\hfil&&\null\hfil&&any smooth\hfil&&any
smooth\hfil&\cr
height4pt&\omit&&\omit&&\omit&&\omit&&\omit&&\omit&&\omit&\cr
&3\hfil&&4\hfil&&2\hfil&&0\hfil&&889\hfil&&any smooth\hfil&&any
smooth\hfil&\cr
height4pt&\omit&&\omit&&\omit&&\omit&&\omit&&\omit&&\omit&\cr
&3\hfil&&\ $\geq 5$\hfil&&1\hfil&&0\hfil&&\null\hfil&&non
trigonal\hfil&&any smooth\hfil&\cr
height4pt&\omit&&\omit&&\omit&&\omit&&\omit&&\omit&&\omit&\cr
&4\hfil&&2\hfil&&1\hfil&&0\hfil&&\null\hfil&&smooth ramif.
div.\hfil&&any smooth\hfil&\cr
height4pt&\omit&&\omit&&\omit&&\omit&&\omit&&\omit&&\omit&\cr
&4\hfil&&3\hfil&&2\hfil&&0\hfil&&1209\hfil&&any smooth\hfil&&any
smooth\hfil&\cr
height4pt&\omit&&\omit&&\omit&&\omit&&\omit&&\omit&&\omit&\cr
&4\hfil&&\ $\geq 4$\hfil&&1\hfil&&0\hfil&&\null\hfil&&any
smooth\hfil&&any
smooth\hfil&\cr
height4pt&\omit&&\omit&&\omit&&\omit&&\omit&&\omit&&\omit&\cr
&\ $\geq 5$\hfil&&2\hfil&&1\hfil&&0\hfil&&\null\hfil&&smooth
ramif.
div.\hfil&&any smooth\hfil&\cr
height4pt&\omit&&\omit&&\omit&&\omit&&\omit&&\omit&&\omit&\cr
&\ $\geq 5$\hfil&&\ $\geq 3$\hfil&&1\hfil&&0\hfil&&\null\hfil&&any
smooth\hfil&&any smooth\hfil&\cr } \hrule}}

%\vskip .5cm

Note that the empty spaces in column 5 of the above table
correspond to the cases were
there are no threefolds with the given $r$ and $g$.

We now prove the main result of the section, which is the
announced computation of the corank of the Gaussian map.

{\bf Theorem (2.13).} {\sl For $r \geq 1$ and $g \geq 2$ let
$S_{r,g}$ be a
smooth K3 surface representing a point of ${\cal H}_{r,g}$ and
$C_{r,g}$ any
smooth hyperplane section of $S_{r,g}$. Then under the hypotheses
on $S_{r,g}$
and $C_{r,g}$ given in columns 6 and 7 of Table (2.12), the values
of
$corank \ \Phi_{\omega_{C_{r,g}}}$ are given in column 3.}

\noindent {\it Proof:} For $r = 1, g \geq 6$ and $C_{r,g}$ general
the
corank of the Gaussian map $\Phi_{\omega_{C_{r,g}}}$ was already
computed in
previous articles: for $6 \leq g \leq 9$ or $g = 11$ it was done
in [CM1]; for
$g = 10$ in [CU]; for $g \geq 12$ in Theorem 4 and Proposition 3
of [CLM1].
Now consider $r$ and $g$ as in the table below
\vskip .3cm
\centerline{Table (2.14)}

\centerline{\vbox{\offinterlineskip
\hrule
\halign{&\vrule#&
\strut\hfil#\tabskip=1.5em\cr
height4pt&\omit&&\omit&\cr
&\ \ $r$\hfil&&\ \ $g$\hfil&\cr
height4pt&\omit&&\omit&\cr
\noalign{\hrule}
height4pt&\omit&&\omit&\cr
&\ \ \ 1\hfil&&\ \ \ \ \ \ $\geq 17$\hfil&\cr
height4pt&\omit&&\omit&\cr
&\ \ \ 2\hfil&&\ \ \ \ \ \ \ $\geq 7$\hfil &\cr
height4pt&\omit&&\omit&\cr
&\ \ \ 3\hfil&&$\geq 5$&\cr
height4pt&\omit&&\omit&\cr
&\ \ \ 4\hfil&&2 or $\geq 4$&\cr
height4pt&\omit&&\omit&\cr
&\ $\geq 5$ \hfil&&$\geq 2$&\cr  }
\hrule}}

\noindent and $C_{r,g}$ any smooth hyperplane section of $S_{r,g}$.
In all
these cases we will prove that $corank \ \Phi_{\omega_{C_{r,g}}} =
1$. By
Theorem (1.1) and Remark (2.7) of [CLM2] we have that the Gaussian
map
$\Phi_{{\cal O}_{S_{r,g}}(1)}$ is surjective for the values of $r$
and $g$ in
Table (2.14). Under the same hypotheses from Lemma (2.3) we have
that
$\phi_{r,g} = \phi_{r,g,1}$ is surjective. Hence Lemma (2.2) gives
that
$\Phi_{\omega_{C_{r,g}}}$ has corank one.

It remains to compute the corank of $\Phi_{\omega_{C_{r,g}}}$ in the
lower
genera. For $g \geq 3$ we have that $S_{r,g} = v_{r}(S_{g})$, where
$v_{r}$ is
the $r$-th Veronese embedding, hence $corank \
\Phi_{\omega_{C_{r,g}}}
= corank \ \Phi_{\omega_{C'}}$ where $C'$ is a smooth curve cut
out on $S_{g}$
by a hypersurface of degree $r$. When $3 \leq g \leq 5$ the K3
surface $S_{g}$
is a complete intersection, hence so is $C'$ and the corank of
the Gaussian map
of $C'$ can be computed by a result of Wahl as follows. If $C'
\subset \I^g$ is
a complete intersection of type $(d_{1}, \ldots , d_{g-1})$,
setting $k =
\sum\limits_{j=1}^{g-1} d_{j} - g - 1$ we have a diagram
%\eject
$$\matrix{\bigwedge^2
H^0(\I^g,{\cal O}_{\I^g}(k))  &
\mapright{\Phi_{{\cal O}_{\I^g}(k)}} &
H^0(\I^g,\Omega_{\I^g}^1(2k)) \cr   && \downarrow \phi_{k} \cr
\downarrow
\pi_{k} & & \ \ H^0(C',\Omega_{\I^g}^1(2k)_{|C'}) \cr
&&\downarrow \psi_{k}
\cr \bigwedge^2 H^0(C',\omega_{C'})
&\mapright{\Phi_{\omega_{C'}}}  &
H^0(C',\omega_{C'}^{\otimes 3}) \cr} \leqno (2.15)$$
\noindent  and $\pi_{k}$ is surjective since
$H^1({\cal I}_{C'}(k)) = 0$,
$\Phi_{{\cal O}_{\I^g}(k)}$ is surjective by [W, Theorem 6.4],
$\phi_{k}$ is
surjective if $2k \not= d_{j}, \ \forall j = 1, \ldots ,
g-1$ ([W, Proposition
6.6]). In the cases at hand the $(g-1)$-tuples $(d_{1}, \ldots ,
d_{g-1})$ are
given by the following

\centerline{\bf Table (2.16)}

\centerline{\vbox{\offinterlineskip
\hrule
\halign{&\vrule#&
\strut\hfil#\tabskip=2em\cr height4pt&\omit&&\omit&&\omit&&\omit&\cr
&\ \
$(d_{1}, \ldots , d_{g-1})$\hfil&&\ \ \ \ $g$\hfil&&\ \ \ \ \
$r$\hfil&&\ \ \ \
$k$\hfil&\cr height4pt&\omit&&\omit&&\omit&&\omit&\cr
\noalign{\hrule} height4pt&\omit&&\omit&&\omit&&\omit&\cr
&$(2,2,2,2)$\hfil&&5&&2&&2&\cr height4pt&\omit&&\omit&&\omit&&\omit&\cr
&$(2,2,3)$\hfil&&4&&2&&2&\cr height4pt&\omit&&\omit&&\omit&&\omit&\cr
&$(2,3,3)$\hfil&&4&&3&&3&\cr height4pt&\omit&&\omit&&\omit&&\omit&\cr
&$(2,4)$\hfil&&3&&2&&2&\cr height4pt&\omit&&\omit&&\omit&&\omit&\cr
&$(3,4)$\hfil&&3&&3&&3&\cr height4pt&\omit&&\omit&&\omit&&\omit&\cr
&$(4,4)$\hfil&&3&&4&&4&\cr } \hrule}}

\noindent thus we see that $2k \not= d_{j}$ (and  $\phi_{k}$ is
surjective)
except for the case $(2,4)$. For the map $\psi_{k} :
H^0({\Omega_{\I^g}^1}_{|C'}(2k)) \to H^0(\Omega_{C'}^1(2k))$ we get
$corank
\ \psi_{k} = H^1(N_{C'/\I^{g}}^*(2k))$ as long as
$H^1({\Omega_{\I^g}^1}_{|C'}(2k)) = 0$. From the Euler sequence
$$0 \to {\Omega_{\I^g}^1}_{|C'}(2k) \to H^0({\cal O}_{C'}(1)) \otimes
{\cal
O}_{C'}(2k-1) \to {\cal O}_{C'}(2k) \to 0$$
\noindent and the fact that $H^1({\cal O}_{C'}(2k-1)) =
H^0({\cal O}_{C'}(1-k))^* = 0$ we see that
$H^1({\Omega_{\I^g}^1}_{|C'}(2k)) = 0$ since the multiplication map
$H^0({\cal O}_{C'}(2k-1)) \otimes  H^0({\cal O}_{C'}(1)) \to H^0({\cal
O}_{C'}(2k))$ is surjective. Excluding the case $(2,4)$ for the moment,
by
(2.15) we have
$$corank \ \Phi_{\omega_{C'}} = corank \ \psi_{k} =
h^1(N_{C'/\I^{g}}^*(2k)) =
h^0(N_{C'/\I^{g}}(-k)) = \sum\limits_{j=1}^{g-1}
h^0({\cal O}_{C'}(d_{j}-k))$$
\noindent and this gives the required values. In the case $(2,4)$,
$C' = Q \cap
S_{3} \subset \I^3$ with $Q$ a quadric surface, and

\noindent {\it Claim} (2.17). {\sl In diagram (2.15) we have:

(i) $corank \ \phi_{2} = 1$;
\eject
(ii) $corank \ \psi_{2} = 10$;

(iii) $Ker \ \psi_{2} \not\subseteq Im \ \phi_{2}$.}

\noindent By Claim (2.17) we deduce that
$corank \ \Phi_{\omega_{C'}} = corank \ \Phi_{\omega_{C'}} \circ
\pi_{2} =
corank \ \psi_{2} \circ \phi_{2} \circ \Phi_{{\cal O}_{\I^3}(2)} =
corank \ \psi_{2} \circ \phi_{2} = corank \ \psi_{2} = 10$.

\noindent {\it Proof of Claim (2.17):} By the above computation
$corank \ \psi_{2} = h^0({\cal O}_{C'}) + h^0({\cal O}_{C'}(2)) =
1 + g(C') =
10$, hence (ii). From the exact sequence
$$0 \to \Omega_{\I^3}^1 \otimes {\cal I}_{C'}(4) \to
\Omega_{\I^3}^1(4)
\to \Omega_{\I^3}^1(4)_{|C'} \to 0$$
\noindent we see that $corank \ \phi_{2} =
h^1(\Omega_{\I^3}^1 \otimes {\cal I}_{C'}(4))$ since
$H^1(\Omega_{\I^3}^1(4)) = 0$. Tensoring by $\Omega_{\I^3}^1(4)$
the Koszul
resolution of ${\cal I}_{C'}$ we have
$$0 \to \Omega_{\I^3}^1(-2) \to \Omega_{\I^3}^1(2) \oplus
\Omega_{\I^3}^1 \to
\Omega_{\I^3}^1 \otimes {\cal I}_{C'}(4) \to 0$$
\noindent and therefore $h^1(\Omega_{\I^3}^1 \otimes {\cal
I}_{C'}(4))
= 1$ by Bott vanishing, hence (i). To see (iii)
consider the following diagram
$$\matrix{0 \to  H^0(N_{S_{3}}^*(4)) \ \ \mapright{f}
\ \ H^0({\Omega_{\I^3}^1}_{|S_{3}}(4)) \cr
\ \ \ \ \ \ \ \ \ \ \ \ \ \ \ \ \ \ \ \ \ \ \  \ \ \ \ \ \ \ \ \
\ \ \ \ \ \ \
\ \ \ \nwarrow{\phi_{3}} \cr \ \ \ \ \ \ \ \ \ \ \ \ \ \ \ \ \ \
\ \ \ \ \ \ \
\ \ \ \ \ \ \ \ \ \ \ \ \ \ \ \ \ \ \ \ \ \ \downarrow{\pi} \ \ \
\ \ \
H^0({\Omega_{\I^3}^1}(4)) \cr
\ \ \ \ \ \ \ \ \ \ \ \ \ \ \ \ \ \ \ \ \ \ \  \ \ \ \ \ \ \ \ \
\ \ \ \ \ \ \
\ \ \ \swarrow{\phi_{2}} \cr   0 \to  H^0(N_{C'}^*(4)) \ \
\mapright{g}
\ \ H^0({\Omega_{\I^3}^1}_{|C'}(4)) \cr}$$

\noindent where $Im \ f = Ker \ \{ \psi_{3} :
H^0({\Omega_{\I^3}^1}_{|S_{3}}(4))
\to H^0({\Omega_{S_{3}}^1}(4)) \} \cong \hbox{\C}$ and $Im \ g =
Ker \
\psi_{2}$ and $\pi$ is the restriction map. By [CLM2, proof of
Theorem (1.1)] we
know that $corank \ \phi_{3} = 1$ and that, if $F$ is a generator of
$Ker \
\psi_{3}$, we have that $f(F) \not\in Im \ \phi_{3}$. If we let
$F' \in H^0(N_{C'}^*(4))$ be such that $g(F') \not= 0$ spans $Ker \
\psi_{2}$
and such that $g(F') = \pi(f(F))$ it follows
that $g(F') \not\in Im \ \phi_{2}$ and hence that  $Ker \ \psi_{2}
\not\subseteq Im \ \phi_{2}$. For otherwise, since  $corank \
\phi_{3} = 1$, we
have that $Im \ \pi \subseteq Im \ \phi_{2}$ and this contradicts
(i) since
$\pi$ is surjective, as it can be easily seen by Bott vanishing.
This proves
Claim (2.17).

In the case $g = 6, r = 2$ and $S_{2,6}$ non trigonal, that is
when
$C' = S_{6} \cap Q \subset \I^6$, with $Q$ a quadric hypersurface
and
$S_{6}$ non trigonal, we have a diagram
%\eject
$$\matrix{\bigwedge^2 H^0(S_{6},{\cal
O}_{S_{6}}(2))  & \mapright{\Phi_{{\cal O}_{S_{6}}(2)}} &
H^0(S_{6},\Omega_{S_{6}}^1(4)) \cr   && \downarrow \phi_{S_{6}} \cr
\downarrow
\pi_{S_{6}} & & \ \ H^0(C',\Omega_{S_{6}}^1(4)_{|C'}) \cr
&&\downarrow
\psi_{S_{6}} \cr \bigwedge^2 H^0(C',\omega_{C'})
&\mapright{\Phi_{\omega_{C'}}}  & H^0(C',\omega_{C'}^{\otimes 3})
\cr}$$
\noindent and

\noindent {\it Claim} (2.18). {\sl In the above diagram we have:

(i) $\pi_{S_{6}}$ and $\Phi_{{\cal O}_{S_{6}}(2)}$ are surjective;

(ii) $corank \ \psi_{S_{6}} \leq 1$;

(iii) $corank \ \phi_{S_{6}} = 1$.}

\noindent From Claim (2.18) it follows that
$corank \ \Phi_{\omega_{C_{2,6}}} = corank \ \Phi_{\omega_{C'}} =
corank \
\psi_{S_{6}} \circ \phi_{S_{6}} \leq 2$. On the other hand if
$corank \ \Phi_{\omega_{C_{2,6}}} \leq 1$ then $H^0(N_{C_{2,6}}(-2)) =
0$ since
$C_{2,6}$ is a canonical curve: this follows by [CM1, Lemma 1.10] in
the case
$corank \ \Phi_{\omega_{C_{2,6}}} = 0$ and by [CM2, Lemma 5.2] in the
case
$corank \ \Phi_{\omega_{C_{2,6}}} = 1$. By a theorem of Zak (see [Z],
[B], [BEL] or [L])
we have that $C_{2,6} \subset \I^{20}$ is not 2-extendable, that is
there is
no threefold $V \subset \I^{22}$ different from a cone such that
$C_{2,6} = V
\cap \I^{20}$. But this contradicts the fact that the K3 surface
$S_{6} \subset
\I^{6}$ is a quadric section of a threefold $V_{3} \subset \I^{6}$ :
$V_{3} =
\hbox{\ms G} \cap H_{1} \cap H_{2} \cap H_{3}$ where $\hbox{\ms G} =
\hbox{\ms
G}(1,4) \subset \I^{9}$ is the Grassmannian of lines in $\I^4$ in the
Pl\"ucker
embedding and the $H_{i}$'s are hyperplanes. Therefore $C_{2,6} =
v_{2}(S_{6}) \cap H' = v_{2}(V_{3}) \cap H' \cap H''$ (with $H', H''$
hyperplanes).

\noindent {\it Proof of Claim (2.18):} As above we have $dim \ Coker
\ \psi_{S_{6}} \leq h^1(N_{C'/S_{6}}^*(4)) = h^1({\cal O}_{C'}(2)) =
h^1(\omega_{C'}) = 1$, that is (ii). For (iii) observe that $dim \
Coker \
\phi_{S_{6}} = h^1(\Omega_{S_{6}}^1(2)) = 1$ by Lemma
(2.3). Also the map $\pi_{S_{6}}$ is surjective since it is the
restriction map
and $H^1({\cal I}_{C'/S_{6}}(2)) = H^1({\cal O}_{S_{6}}) = 0$. The
fact that
$\Phi_{{\cal O}_{S_{6}}(2)}$ or equivalently
$\Phi_{{\cal O}_{S_{2,6}}(1)}$ is
surjective has been proved in [CLM2, proof of Theorem (1.1)] and
this gives (i).
This ends the proof of the Claim (2.18).

To finish the proof of the theorem let us consider the cases $g = 2,
r = 2, 3$.
For $g = 2$ the K3 surface $S_{2}$ is a double cover of $\I^2$ given
by a
linear system $|H|$ with $H^2 = 2$ and the curve $C' \cong C_{r,2}$ is
a smooth
member of $|rH|$. For $r = 2$ we have that $C_{2,2}$ is a hyperelliptic
curve
of genus $5$, hence  $corank \ \Phi_{\omega_{C_{2,2}}} = 13$ by [CM2,
Proposition 1.1] or [W, Remark 5.8.1]. When $r = 3$ a general member
$C_{3,2}$
of $|3H|$ is isomorphic to a smooth plane sextic, hence $corank \
\Phi_{\omega_{C_{3,2}}} = 10$ by [CM2, section 2], while the special
members
have a 2:1 map onto a smooth plane cubic, hence  $corank \
\Phi_{\omega_{C_{3,2}}} = 18$ by [CM2, Corollary 3.3].

\noindent This then concludes the proof of Theorem (2.13). \ \qed

{\bf (2.19)} {\it Remark.} Note that in the case $r = 3, g = 2$ of
Theorem (1.1), the corank of the Gaussian map $\Phi_{\omega_{C_{3,2}}}$
{\it is
not constant} in the linear system $|{\cal O}_{S_{3,2}}(1)|$. On the
other hand
using a diagram like (2.1) it is easy to deduce, as in Lemma (2.2),
that if $S$
is a smooth K3 surface, $L$ is an effective line bundle on $S$ and
we have
$\Phi_{L}$ surjective and $H^1(S,\Omega_{S}^1 \otimes L) = 0$, then,
for every
smooth $C \in |L|$, the corank of $\Phi_{\omega_{C}}$ is one. It is an
interesting question to find better hypotheses that insure
the constancy of the corank of $\Phi_{\omega_{C}}$ in a given linear
system on a
smooth K3 surface. The other intriguing feature is that the example
given
above of non constancy of the corank of the Gaussian map is the same
as
Donagi-Morrison's example of non constancy of the gonality of the
smooth
curves in a linear system on a K3 surface.
\vskip .3cm
Before the end of this section we will also prove a result about the
Gaussian
maps $\Phi_{\omega_{C_{r,g}},\omega_{C_{r,g}}^{\otimes k}}$ that will
be useful
in the next section.

{\bf Proposition (2.20).} {\sl Fix $r \geq 1$ and $g \geq 2$. Assume
that
$S_{r,g}$ and $C_{r,g}$ satisfy the hypotheses given in columns 6 and
7 of
Table (2.12). Then the Gaussian map
$\Phi_{\omega_{C_{r,g}},\omega_{C_{r,g}}^{\otimes k}}$ is surjective
for every
$k \geq 2$ except for $k = 2$ and $(r,g) = (2,3), (1,6), (3,2)$.
Moreover if
$N_{C_{r,g}}$ is the normal bundle of $C_{r,g}$ in $\I H^0(C_{r,g},
{\cal
O}_{C_{r,g}}(1))^*$, then $h^0(N_{C_{r,g}}(-2))$ is given by column 4
of Table
(2.12) and $h^0(N_{C_{r,g}}(-k)) = 0$ for every $k \geq 3$.}

\noindent {\it Proof:} Set $C = C_{r,g}$ and $S = S_{r,g}$. Since $C$ is
a
canonical curve it is well known that, for $k \geq 2$, we have
$$Coker \ \Phi_{\omega_{C},\omega_{C}^{\otimes k}} \cong
H^0(N_{C}(-k))^*$$
\noindent (see for example [CM1, Proposition 1.2]), hence the assertion
on
surjectivity of $\Phi_{\omega_{C_{r,g}},\omega_{C_{r,g}}^{\otimes k}}$
follows by proving the vanishing of $h^0(N_{C_{r,g}}(-k))$. For $r = 1,
g
\geq 6$ and $C$ general the values of $h^0(N_{C}(-k))$ are given in
Lemma 4 of
[CLM1]. Suppose now $r = 1, g \geq 17$ or $r \geq 2, g \geq 3$ or $r
\geq 5, g
= 2$ and $C$ any smooth hyperplane section of $S$. Set $R({\cal
O}_{S}(1),{\cal O}_{S}(k)) = Ker \ \{ H^0({\cal O}_{S}(1)) \otimes
 H^0({\cal
O}_{S}(k)) \to H^0({\cal O}_{S}(k+1)) \}$ and consider the diagram
%\eject
$$\matrix{R({\cal O}_{S}(1),{\cal
O}_{S}(k))  & \mapright{\Phi_{{\cal O}_{S}(1),{\cal O}_{S}(k)}} &
H^0(S,\Omega_{S}^1(k+1)) \cr   && \downarrow \phi \cr  \downarrow \pi
& & \ \
H^0(C,\Omega_{S}^1(k+1)_{|C}) \cr  &&\downarrow \psi \cr  R({\cal
O}_{C}(1),{\cal O}_{C}(k))
&\mapright{\Phi_{\omega_{C},\omega_{C}^{\otimes
k}}}  &  H^0(C,\omega_{C}(k+1)). \cr}$$
\noindent Since $\Phi_{{\cal O}_{S}(1),{\cal O}_{S}(k)}$ is surjective
by
[CLM2, Theorem (1.1) and Remark (2.7)], and $\psi$ is surjective because
$Coker
\ \psi \subseteq H^1(N_{C/S}^*(k+1)) = H^1({\cal O}_{C}(k)) = 0$ for $k
\geq 2$,
it follows that $$corank \ \Phi_{\omega_{C},\omega_{C}^{\otimes k}}
\leq corank
\ \phi =  h^1(S,\Omega_{S}^1(k)) =
\cases{ $0$ &for $(k,r,g) \not= (2,2,3)$  \cr
$1$ &for $(k,r,g) = (2,2,3)$ \cr}$$
\noindent by Lemma (2.3).

For $g = 2, r = 3$ the values of $h^0(N_{C}(-k))$ are in [CM2, section
2] for
$C$ general and in [CM2, Theorem (3.2)] for $C$ special. For $g = 2,
 r = 4$ we
have $corank \ \Phi_{\omega_{C}} = 1$ by Theorem (2.13), hence
$h^0(N_{C}(-k)) = 0$ for $k \geq 2$ ($k = 2$ by [CM2, Lemma 5.2]).
\ \qed

{\bf (2.21)} {\it Remark.} In the case $r = 1, g \geq 17$ Theorem
(2.13) is a
generalization of the Corank One Theorem of [CLM1] (Theorem 4)
since there the
fact that  $corank \ \Phi_{\omega_{C_{g}}} = 1$ was proved for a
{\it general}
hyperplane section $C_{g}$ of a {\it general} prime K3 surface
$S_{g}$, while
here the same holds for {\it any} smooth hyperplane section $C_{1,g}$
of a
general $S_{1,g}$. On the other hand it should be noticed that the
proof of
this fact is not independent of Theorem 4 of [CLM1] as the latter is
used to
deduce the Corank One Theorem from the surjectivity of $\Phi_{{\cal
O}_{S_{1,g}}(1)}$ (in Lemma (2.3)). Alternatively one can avoid using
Theorem 4
of [CLM1] by invoking results of Mori and Mukai ([MM], [M1]); see
remark (2.11).

\vskip .5cm

\noindent {\bf 3. VARIETIES WITH CANONICAL CURVE SECTION}
\vskip .5cm

We begin by recalling Zak's theorem mentioned in the introduction.
A smooth
nondegenerate variety $X \subset \I^m$ is said to be $k${\it-extendable}
if
there exists a variety $Y \subset \I^{m + k}$ that is not a cone and
such that
$X = Y \cap \I^m$. Zak's theorem says that if $codim \ X \geq 2$,
$h^0(N_{X}(-1)) \leq m + k$ and $h^0(N_{X}(-2)) = 0$, then $X$ is not
k-extendable. For a canonical curve $C \subset \I^{g-1}$ one has also
$h^0(N_{C}(-1)) = g + \hbox{corank} \ \Phi_{\omega_{C}}$ ([CM1,
Proposition
1.2]), hence the knowledge of $corank \ \Phi_{\omega_{C}}$ (and of
$h^0(N_{C}(-2))$ that is also given by a Gaussian map) gives
information on the
possibility of extending $C$.

Next we consider a smooth nondegenerate variety $X \subset \I^m$ of
dimension
$n \geq 3$ such that its general curve sections $C$ are canonical
curves. By a
well-known equivalence criterion we must have $- K_{X} = (n-2)H$,
where $H$ is the hyperplane divisor. Hence a general surface section
$S$ is a
K3 surface, and therefore $S = S_{r,g}, C = C_{r,g}$ for some $r,g$.
We denote
such an $X$ by $X_{r,g}^n$. When $n = 3$ we denote $V_{r,g} =
X_{r,g}^3$. Note
that $V_{r,g} \subset \I^{N(r,g) + 1} =  \I H^0(V_{r,g},-
K_{V_{r,g}})^*$ is
anticanonically embedded and, by Riemann-Roch, $N(r,g) =
h^0(V_{r,g},-
K_{V_{r,g}}) - 2 = - {1 \over 2} K_{V_{r,g}}^3 + 1 = 1 + r^2(g-1) =
{1 \over
2}d r^3 + 1$, where $H = r \Delta, d = \Delta^3$. Hence to study
varieties with
canonical curve sections it is not restrictive to study the
varieties
$X_{r,g}^n$.

These varieties $X_{r,g}^n$ have been extensively studied, both in
dimension three (Fano, Iskovskih [I1], [I2], etc.) and in dimension
$n \geq
4$ (Mukai [M1], [M2]). Observe that by applying some basic adjunction
theory
(e.g. [KO]) we can divide the varieties $X_{r,g}^n$ in two categories
(with two
exceptions): (a) $n = 3$ and $r \geq 1$, that is Fano threefolds of
index $r$;
(b) $n \geq 4$ and $r = 1$, that is Fano manifolds of dimension $n
\geq 4$ and
index $1$. In fact for $n \geq 4, r \geq 2$ we have only $n = 4, 5,
r = 2$ and
$(X, \Delta) = (Q, {\cal O}_{Q}(1)), (\I^5, {\cal O}_{\I^5}(1))$,
where $Q$ is
a smooth quadric in $\I^5$, since in all other cases it is $r(n-2) >
n+1$. Even
though we will not make use of adjunction theory in our proofs, the
interesting
cases are (a) and (b) above.

The case of Fano threefolds of index one was already treated in
[CLM1], thus from now on we will assume $r \geq 2$ or $n \geq 4$.

Given $r \geq 2$ and $g \geq 2$, we will denote by ${\cal V}_{r,g}$
the Hilbert scheme of smooth Fano threefolds $V_{r,g} \subset
\I^{N(r,g) + 1}$
of index $r$, anticanonically embedded and by ${\cal V}_{r,g,1}$
the closure in
${\cal V}_{r,g}$ of the locus of Fano threefolds with Picard number
one.

We now present a list of examples taken from the papers of Fano and
Iskovskih of
families of Fano threefolds $V_{r,g}$ of the principal series and of
index
greater than one (see e.g. [M]):

\centerline{\bf Table (3.1)}

\centerline{\vbox{\offinterlineskip
\hrule
\halign{&\vrule#&
\strut\hfil#\tabskip=.05em\cr
height4pt&\omit&&\omit&&\omit&&\omit&&\omit&&\omit&\cr
&\ $r$\ &&\ $g$ &&\ $d = \Delta^3$\hfil&& \ $N(r,g)$&& \ number
of
parameters\hfil&&$V_{r,g}$ \hfil&\cr
height4pt&\omit&&\omit&&\omit&&\omit&&\omit&&\omit&\cr
\noalign{\hrule}
height4pt&\omit&&\omit&&\omit&&\omit&&\omit&&\omit&\cr
&4\hfil&&3\hfil&&1\hfil&&33\hfil&&1209\hfil&&$\I^3$\hfil&\cr
height4pt&\omit&&\omit&&\omit&&\omit&&\omit&&\omit&\cr
&3\hfil&&4\hfil&&2\hfil&&28\hfil&&889\hfil&&quadric in
$\I^4$\hfil&\cr
height4pt&\omit&&\omit&&\omit&&\omit&&\omit&&\omit&\cr
&2\hfil&&3\hfil&&2\hfil&&9\hfil&&139\hfil&&\ double cover of
$\I^3$ ramified
along a
quartic \hfil&\cr
height4pt&\omit&&\omit&&\omit&&\omit&&\omit&&\omit&\cr
&2\hfil&&4\hfil&&3\hfil&&13\hfil&&234\hfil&&cubic in
$\I^4$\hfil&\cr
height4pt&\omit&&\omit&&\omit&&\omit&&\omit&&\omit&\cr
&2\hfil&&5\hfil&&4\hfil&&17\hfil&&363\hfil&&\ complete
intersection of two quadrics in
$\I^5$\hfil&\cr
height4pt&\omit&&\omit&&\omit&&\omit&&\omit&&\omit&\cr
&2\hfil&&6\hfil&&5\hfil&&21\hfil&&525\hfil&&$\hbox{\ms G}(1,4)
\cap \I^6 \subset
\I^9$\hfil&\cr }  \hrule}}

Our first classification result is the following:
\vskip .3cm
{\bf Theorem (3.2).} {\sl Fix $r \geq 2$ and $g \geq 2$. Then we
have:

(i) ${\cal V}_{r,g,1} = \emptyset$ for $(r,g)$ not in Table (3.1);

(ii) ${\cal V}_{r,g} = \emptyset$ for $(r,g)$ such that: $r = g =
2; r = 3,
g = 2, 3; r = 4, g = 2$ or $g \geq 4; r \geq \indent 5, g \geq 2;
r = 3, g \geq
5$ assuming furthermore that the hyperplane section $S_{3,g}$ of
$V_{3,g}$ is
\indent non trigonal;

(iii) For $(r,g)$ as in Table (3.1) ${\cal V}_{r,g,1}$ is
irreducible and the
families of Fano-Iskovskih \indent form a dense open subset of
smooth points
of ${\cal V}_{r,g,1}$.}

\noindent {\it Proof:} Suppose first $g = 2$. Then $d = 2(g-1)/r =
2/r$ is an
integer only if $r = 2$, but then for the K3 surface $S_{2,2}$ we
have
that  ${\cal O}_{S_{2,2}}(1) = {\cal O}_{S_{2,2}}(2 \Delta_{|S_{2,2}})
\cong
{\cal O}_{S_{2}}(2D)$ where $\phi_{D}$ gives the 2:1 map onto $\I^2$.
But this
is a contradiction since on $S_{2}$ the line bundle
${\cal O}_{S_{2}}(2D)$ is not very
ample. Similarly for $r = g = 3$ we have that $d$ is not an integer.
For $r = 4, g
\geq 4$ or $r \geq 5, g \geq 3$ or $r = 3, g \geq 5$ and $S_{3,g}$
non trigonal, we
have by Theorem (2.13) that $corank \ \Phi_{\omega_{C_{r,g}}} = 1$,
hence by [CM1,
Proposition 1.2]  $h^0(N_{C_{r,g}}(-1)) = N(r,g) + 1$. Also
$h^0(N_{C_{r,g}}(-2)) = 0$
by [CM2, Lemma 5.2], and therefore Zak's theorem implies that
$C_{r,g}$ is not
$2$-extendable, hence that ${\cal V}_{r,g} = \emptyset$ for the
above values of
$r$ and $g$. This proves (ii). To see (i) and (iii) we will use
the following

{\bf Lemma (3.3).} {\sl Suppose $V_{r,g}$ has Picard number one.
Then
a general hyperplane section $S_{r,g}$ of a general deformation of
$V_{r,g}$ is
represented by a general point of ${\cal H}_{r,g}$.}

\noindent {\it Proof of Lemma (3.3):} Set $V = V_{r,g}, S = S_{r,g}$
and
consider the usual deformation diagram
$$\matrix{H^1(V,T_{V})  & \mapright{\alpha_{S}} & H^1(V,T_{V_{|S}})
& \to
H^2(V,T_{V}(-1)) & \mapright{\gamma_{S}} & H^2(V,T_{V}) & \to
H^2(V,T_{V_{|S}}) \cr
&& \uparrow \beta_{S} \cr  && H^1(S,T_{S}) \cr} $$
\noindent First we prove

(3.4) $\beta_{S}$ is surjective;

(3.5) $h^2(V,T_{V}(-1)) = 1$;

(3.6) $H^2(V,T_{V}) = 0$.

\noindent In  fact the map $\beta_{S}$ comes from the exact
sequence
$$0 \to T_{S} \to T_{V_{|S}} \to N_{S/V} \to 0 \leqno (3.7)$$
\noindent hence $Coker \ \beta_{S} \subseteq H^1(N_{S/V}) =
H^1({\cal
O}_{S}(1)) = 0$, that is (3.4). By Serre duality we get
$h^2(V,T_{V}(-1)) = h^1(\Omega_{V}^1) = 1$ since ${\cal O}_{V}(1)
= {\cal
O}_{V}(-K_{V})$ and $V$ has Picard number one. This proves (3.5).
To see
(3.6) observe that by (3.7) we have $H^2(T_{V{|S}}) = 0$, hence if
$H^2(T_{V}) \not= 0$ we necessarily have $h^2(T_{V}) = 1$ and
$\gamma_{S}$ is
surjective. On the other hand $\gamma_{S}$ is given by multiplication
with the
equation of the hyperplane defining $S$, hence the multiplication map
$\mu :
H^0({\cal O}_{V}(1)) \otimes H^2(T_{V}(-1)) \to H^2(T_{V})$ is
surjective. But
then, by Bertini's theorem, there is a smooth surface $S'$ in the
codimension one sublinear system $Ker \ \mu \subset
H^0({\cal O}_{V}(1))$.
Replacing $S$ with $S'$ we get that $\gamma_{S'}$ is the zero
map, and this
contradicts the hypothesis that $H^2(T_{V}) \not= 0$.
Alternatively (3.6)
follows by Kodaira-Nakano vanishing since $H^2(T_{V}) \cong
H^1(\Omega_{V}^1
(K_{V})) = H^1(\Omega_{V}^1(-1)) = 0$.

Now by (3.6) every infinitesimal deformation of $V$ is
unobstructed and by
(3.5) we have that $dim \ Coker \ \alpha_{S} = 1$, hence
every infinitesimal
deformation of $S$ comes from a deformation of $V$. This proves
Lemma (3.3).

To see (i) of Theorem (3.2) by Lemma (3.3) we can assume that
$S_{r,g}$ is
general. Moreover note that if $(r,g)$ are not in Table (3.1) and
not in (ii),
then we have $r = 2, g \geq 7$. Therefore by Theorem (2.13) we have
again
$corank \ \Phi_{\omega_{C_{r,g}}} = 1$, hence we conclude as above
that ${\cal
V}_{r,g,1} = \emptyset$. This gives (i).

To see (iii) observe that the Fano threefolds $V_{r,g}$ of the
principal series
are projectively Cohen-Macaulay, hence they flatly deform to
the cone
$\widehat S_{r,g}$ over their hyperplane section $S_{r,g}$, and
similarly for
$S_{r,g}$ to the cone $\widehat C_{r,g}$ over $C_{r,g}$. Setting
$V = V_{r,g},
S = S_{r,g}, C = C_{r,g}, N = N(r,g)$ and denoting by $T_{r,g}$
the tangent
space to ${\cal V}_{r,g,1}$ at the point representing $V_{r,g}$,
we have
$$dim \ T_{r,g} = h^0(N_{V/ \I^{N+1}}) \leq h^0(N_{{\widehat S} /
\I^{N+1}}) = \sum\limits_{h \geq 0} h^0(N_{S/ \I^N}(-h)) \leq
\leqno (3.8)$$
$$ \leq h^0(N_{S/ \I^N})
+ \sum\limits_{h \geq 1} h^0(N_{{\widehat C}/ \I^N}(-h)) \leq
h^0(N_{S/ \I^N}) + \sum\limits_{h \geq 1, j \geq 0} h^0(N_{C /
\I^{N-1}}(-h-j))
=$$
$$ = h^0(N_{S/ \I^N}) + h^0(N_{C / \I^{N-1}}(-1)) + 2 h^0(N_{C /
\I^{N-1}}(-2))$$
\noindent by Proposition (2.20). As we know $h^0(N_{S/ \I^N}) =
dim \ {\cal H}_{r,g} = 18 + (N+1)^2$ and $h^0(N_{C / \I^{N-1}}(-1))
= N +
corank \ \Phi_{\omega_{C}}$ hence we deduce that
$$h^0(N_{S/ \I^N}) + h^0(N_{C / \I^{N-1}}(-1)) + 2 h^0(N_{C /
\I^{N-1}}(-2)) \leq f(r,g) := \leqno (3.9)$$
$$= N(r,g)^2 + 3 N(r,g) + 19 + corank \
\Phi_{\omega_{C_{r,g}}} + 2 h^0(N_{C_{r,g}}(-2)).$$
\noindent By Theorem (2.13), Proposition (2.20) and Table (3.1) we
see that
$f(r,g)$ is equal to the number of parameters of the examples of
Fano-Iskovskih. By (3.8) and (3.9) we also have that $f(r,g)$ is an
upper
bound for $h^0(N_{{\widehat S} / \I^{N+1}})$, hence the argument
presented in
the Introduction applies and we conclude that ${\cal V}_{r,g,1}$ is
irreducible
and that the families of Fano-Iskovskih form a dense open subset of
smooth
points. \ \qed

We now turn to the higher dimensional case and classify the Mukai
varieties.
For $n \geq 4, r \geq 1$ and $g \geq 2$ we let ${\cal X}_{n,r,g}$ be
the
Hilbert scheme of smooth Fano $n$-folds $X_{r,g}^n \subset
\I^{N(r,g) + n - 2}$
of index $r(n-2)$ and ${\cal X}_{n,r,g,1}$ the closure in
${\cal X}_{n,r,g}$ of
the locus of Fano $n$-folds with Picard number one.

For Fano manifolds of dimension $n \geq 4$ and index $r(n-2)$, of the
principal
series and with Picard number one, we present the list given by Mukai
([M1], [M2]) of
examples of families of such  $X_{r,g}^n$. Of course, since the general
hyperplane section
of an $X_{r,g}^n$ is an $X_{r,g}^{n-1}$, we only present the list of
the maximum
dimensional varieties. There is one family for each genus $g\geq 6$
(for $r = 1, g \leq 5$ they are complete intersections)
and we denote the maximum dimension of the examples in the list by
$n(g)$:

\centerline{\bf Table (3.10)}

\centerline{\vbox{\offinterlineskip
\hrule
\halign{&\vrule#&
\strut\hfil#\tabskip=1em\cr
height4pt&\omit&&\omit&&\omit&&\omit&&\omit&\cr
& \ \ $r$&& \ \ $g$&& \ \ \ $n(g)$&& \ \ number of parameters
\hfil&&$X_{r,g}^{n(g)}$\hfil&\cr
height4pt&\omit&&\omit&&\omit&&\omit&&\omit&\cr
\noalign{\hrule}
height4pt&\omit&&\omit&&\omit&&\omit&&\omit&\cr
&1&&6&&6&&145\hfil&&$\widetilde{\hbox{\ms G}(1,4)} \cap Q
\subset \I^{10}$&\cr
height4pt&\omit&&\omit&&\omit&&\omit&&\omit&\cr
&1&&7&&10&&210\hfil&&\ $SO(10, \hbox{\C}) / P \subset \I^{15}$&\cr
height4pt&\omit&&\omit&&\omit&&\omit&&\omit&\cr
&1&&8&&8&&189\hfil&&$\hbox{{\ms G}(1,5)} \subset \I^{14}$&\cr
height4pt&\omit&&\omit&&\omit&&\omit&&\omit&\cr
&1&&9&&6&&174\hfil&&$Sp(6, \hbox{\C}) / P \subset \I^{13}$&\cr
height4pt&\omit&&\omit&&\omit&&\omit&&\omit&\cr
&1&&10&&5&&181\hfil&&$G_{2} / P \subset \I^{13}$&\cr
height4pt&\omit&&\omit&&\omit&&\omit&&\omit&\cr
&2&&5&&5&&405\hfil&&$v_{2}(\I^5) \subset \I^{20}$&\cr }
\hrule}}

\noindent For $4 \leq n < n(g)$ Mukai shows that any $X_{1,g}^n$ is
a linear
section of the $X_{1,g}^{n(g)}$.

Our final result gives the classification as follows:
\vskip .3cm
{\bf Theorem (3.11).} {\sl Fix $n \geq 4, r \geq 1$ and $g \geq 6$ or
$r \geq 2$
and $g \geq 2$. Then we have:

(i) ${\cal X}_{n,r,g,1} = \emptyset$ for $(r, g)$ not in Table (3.10)
or for
$(r, g)$ in Table (3.10) and $n > n(g)$;

(ii) ${\cal X}_{n,r,g} = \emptyset$ for $(r,g)$ such that: $r = 2,
3$ and $2
\leq g  \leq 4; r \geq 4, g \geq 2$; $S_{r,g}$ non \indent trigonal
and $r = 2, g = 6$,
or $r = 2, g = 5, n \geq 6$, or $r = 3, g \geq 5$;

(iii) For $(r, g)$ as in Table (3.10) and $4 \leq n \leq n(g)$ we
have that
${\cal X}_{n,r,g,1}$ is irreducible \indent and the examples of
families of
Mukai form a dense open subset of smooth points of \indent
${\cal X}_{n,r,g,1}$.}

\noindent {\it Proof:} By (ii) of Theorem (3.2) to prove (ii)
of Theorem (3.11) we need to consider the cases $r = 2$ and $3 \leq g
\leq 6$,
$r = 3$ and $g = 4$ or $r = 4$ and $g = 3$. For $(r,g) = (3,4)$ and
$(r,g)=(4,3)$ we have by Theorem (2.13) that $corank \
\Phi_{\omega_{C_{r,g}}} =
2$ and $h^0(N_{C_{r,g}}(-2)) = 0$ by Proposition (2.20), hence
$C_{r,g}$ is not
$3$-extendable by Zak's theorem. Similarly for $r = 2, g = 5$
(respectively
$g = 6$) and $S_{r,g}$ non trigonal, it is $corank \
\Phi_{\omega_{C_{2,5}}} =
4$ (respectively $corank \ \Phi_{\omega_{C_{2,6}}} = 2$)
and $h^0(N_{C_{r,g}}(-2)) = 0$, hence $C_{2,5}$ is not
$5$-extendable
(respectively $C_{2,6}$ is not $3$-extendable). For $r =
2, g = 4$ we have $deg \ X_{2,4}^n =  deg \ C_{2,4} = 2 N(2,4) -
2 = 24$ and
${\cal O}_{X_{2,4}^n}(1) = {\cal O}_{X_{2,4}^n}(2 \Delta)$, hence
$2^n \Delta^n
= 24$ and therefore $n \leq 3$. For $r = 2, g = 3$ if there
exists an
$X_{2,3}^4$ we get $deg \ X_{2,3}^4 =  deg \ C_{2,3} = 2 N(2,3) - 2 =
16 = 2^4
\Delta^4$, hence $\Delta^4 = 1$. By restricting to the surface
section
$S_{2,3}$ one easily checks that $|\Delta|$ is base point free,
birational and
$h^0({\cal O}_{X_{2,3}^4}(\Delta)) = 4$ and this of course
is a contradiction. Alternatively we can exclude the
existence of a smooth $X_{2,3}^4$ using some basic adjunction
theory: we have
$-K_{X_{2,3}^4} = 4 \Delta$, with $H = 2 \Delta$, hence $\Delta$
is ample; by
[KO] it must be $(X_{2,3}^4, {\cal O}_{X_{2,3}^4}(\Delta) ) \cong
(Q, {\cal
O}_{Q}(1) )$ where $Q$ is a smooth quadric in $\I^5$, but this
contradicts the
fact that $h^0({\cal O}_{X_{2,3}^4}(1)) =
h^0({\cal O}_{X_{2,3}^4}(2\Delta)) =
12$ while $h^0({\cal O}_{Q}(2)) = 20$. This proves (ii).

To see part (i), by part (ii) which was just proved and by (i)
of Theorem
(3.2)(and by Lemma (3.3)), it remains to consider the cases  $r =
1, g \geq 6$.
Let $r = 1, 7 \leq g \leq 10$. We notice that
by [CLM1, Lemma 4] one has $h^0(N_{C_{1,g}}(-2)) = 0$, hence Zak's
theorem
applies, and therefore if we set $\nu:=corank \
\Phi_{\omega_{C_{1,g}}} +1$,
then $C_{1,g}$ is not $\nu$-extendable. On the other hand
$\nu$ is computed in
[CLM1, Table 2] and is equal to $n(g)$. Hence we conclude that
${\cal
X}_{n,r,g,1} = \emptyset$ for $r = 1, 7 \leq g \leq 10$ and $n >
n(g)$.
For $r = 1$ and $ g \geq 11$ by [CLM1, Theorem 4, Lemma 4 and
Table 2] we have
$corank \ \Phi_{\omega_{C_{1,g}}} \leq 2$ and
$h^0(N_{C_{1,g}}(-2)) = 0$, hence
$C_{1,g}$ is not $3$-extendable.

For $r = 1, g = 6$ we have $corank \ \Phi_{\omega_{C_{1,6}}} =
10$ by [CLM1,
Table 2] but in fact there is no smooth $X_{1,6}^7$ in this case
(note that
$h^0(N_{C_{1,6}}(-2)) \not= 0$ !). To see this observe first that,
by (iii)
below, any smooth $X_{1,6}^6$ is the intersection of a quadric
hypersurface $Q
\subset \I^{10}$ and the cone $\widetilde{\hbox{\ms G}}$ over the
Grassmannian
$\hbox{\ms G} = \hbox{\ms G}(1,4) \subset \I^9$ in the Pl\"ucker
embedding. Now
suppose there is $X_{1,6}^7 \subset \I^{11}$ such that $X_{1,6}^6 =
{\widetilde{\hbox{\ms G}}} \cap Q = X_{1,6}^7 \cap H$. Since the
ideal of
$\hbox{\ms G} \subset \I^9$ is generated by five quadrics we have
$H^0({\cal
I}_{\hbox{\ms G}}(2)) = < Q_{1}, \ldots , Q_{5} >$, and from this
it follows
that also the ideals of $X_{1,6}^6$ and $X_{1,6}^7$ are generated
by five
quadrics. In fact if we set $\I^9 = \{x_{10} = x_{11} = 0 \} \subset
\I^{11}$ we
have  $H^0({\cal I}_{X_{1,6}^6}(2)) = < Q_{1},
\ldots , Q_{5}, Q >$ and $H^0({\cal I}_{X_{1,6}^7}(2)) =
< \overline{Q}_{1}, \ldots , \overline{Q}_{5}, \overline{Q} >$,
where
$\overline{Q}_{i}$ and $\overline{Q}$ are quadrics in $\I^{11}$ that
restrict to
$Q_{i}$ and $Q$ respectively for $x_{11} = 0$. But then
the variety $Y \subset \I^{11}$ of dimension $8$ and degree $5$
defined
scheme-theoretically by the five quadrics $\overline{Q}_{i}, i = 1,
\ldots ,
5$, is a variety that extends twice the Grassmannian $\hbox{\ms G}
\subset \I^9$
and hence $Y$ must be a cone over {\ms G} with vertex a line (this
follows for
example from [CLM3]), and then $X_{1,6}^7$ is singular. Therefore
(i) is proved.

As in the proof of part (iii) of Theorem (3.2), using successive
degenerations
to cones over hyperplane sections, to see part (iii) of Theorem
(3.11)
we need to show that the family of known examples of $X_{r,g}^n$
has the same
dimension as the dimension of the tangent space $T_{r,g}^n$ to the
Hilbert scheme
${\cal X}_{n,r,g,1}$ at the cone points. In this case we have
$$dim \ T_{r,g}^n = h^0(N_{X_{r,g}^n}) \leq h^0(N_{V_{r,g}}) +
(n-3) h^0(N_{V_{r,g}}(-1)) + [n - 3 + {n-3 \choose 2}]
h^0(N_{V_{r,g}}(-2))$$
\noindent since $h^0(N_{V_{r,g}}(-k)) = 0$ for $k \geq 3$ because
$h^0(N_{C_{r,g}}(-k)) = 0$ by [CLM1, Lemma 4] and Proposition (2.20).
Therefore
we get
$$h^0(N_{V_{r,g}}) +
(n-3) h^0(N_{V_{r,g}}(-1)) + [n - 3 + {n-3 \choose 2}]
h^0(N_{V_{r,g}}(-2))
\leq $$
$$\leq dim \ {\cal V}_{r,g,1} +  (n-3) h^0(N_{C_{r,g}}(-1)) + [3n -
9 +
{n-3 \choose 2}] h^0(N_{C_{r,g}}(-2)) =$$
$$= N(r,g)^2 + n N(r,g) +
19 + (n - 2) \ corank \ \Phi_{\omega_{C_{r,g}}} + [3n - 7 + {n-3
\choose
2}] h^0(N_{C_{r,g}}(-2)).$$
\noindent Now one easily sees that the values of
the upper bound given above coincide with the number of parameters
of the
examples of Mukai in Table (3.10) for $n = n(g)$ or of their
hyperplane sections
for $n < n(g)$. \ \qed

{\bf (3.12)} {\it Remark.} One interesting consequence of the two
theorems just
proved is that, in general, the converse of Zak's theorem, at least
for {\it
smooth extensions}, does not hold. That is there are smooth
nondegenerate
varieties $X \subset \I^m$ such that $h^0(N_{X}(-2)) = 0$ and $X$
is not
smoothly $k${\it-extendable}, but $h^0(N_{X}(-1)) > m + k$. Examples
are given
by the curves $C_{r,g}$ for $(r,g) = (2,4)$ and $(3,3)$.
Also note that in the case $r = 1, g = 6$ we have that $C_{1,6}$ is
{\it infinitely many times extendable} (just take quadric sections of
cones over
$\hbox{\ms G}(1,4) \subset \I^9$ with vertices linear spaces), but it
is not
smoothly $6${\it-extendable}. Of course in this case we have
$h^0(N_{C_{1,6}}(-2)) \not= 0$.
\vskip .5cm
{\bf REFERENCES}
\baselineskip 12pt
\vskip .5cm
\item{[B]} Badescu,L.:\ On a result of Zak-L'vovsky.\ {\it
Preprint No.\ 13/1993, Institutul de Matematica al Academiei Romane}.
\vskip .2cm
\item{[BEL]} Bertram,A., Ein,L., Lazarsfeld,R.:\ Surjectivity of
Gaussian maps
for line bundles of large degree on curves.\ In:\ {\it Algebraic
Geometry,
Proceedings Chicago 1989.\ Lecture Notes in Math.\ \bf 1479}.\
Springer,
Berlin-New York:\ 1991, 15-25.
\vskip .2cm
\item{[CHM]} Ciliberto,C., Harris,J., Miranda,R.:\ On the surjectivity
of the
Wahl map.\ {\it Duke Math.\ J.\ \bf 57}, (1988) 829-858.
\vskip .2cm
\item{[CLM1]} Ciliberto,C., Lopez,A.F., Miranda,R.:\ Projective
degenerations of K3 surfaces, Gaussian maps and Fano threefolds.\ {\it
Invent.\ Math.\ \bf 114}, (1993) 641-667.
\vskip .2cm
\item{[CLM2]} Ciliberto,C., Lopez,A.F., Miranda,R.:\ On the corank of
Gaussian
maps for general embedded K3 surfaces.\ To appear in {\it Israel
Mathematical Conference Proceedings. Papers in honor of Hirzebruch's
65th
birthday}, vol. \ 9, AMS Publications (1994).
\vskip .2cm
\item{[CLM3]} Ciliberto,C., Lopez,A.F., Miranda,R.:\ Some remarks
on the
obstructedness of cones over curves of low genus.\ {\it In
preparation}.
\vskip .2cm
\item{[CM1]} Ciliberto,C., Miranda,R.:\ On the Gaussian
map for canonical curves of low genus.\ {\it Duke Math.\ J.\ \bf 61},
(1990)
417-443.
\vskip .2cm
\item{[CM2]} Ciliberto,C., Miranda,R.:\ Gaussian maps for certain
families of canonical curves.\ In: {\it Complex Projective Geometry,
Proceedings
Trieste-Bergen 1989.\ London Math.\ Soc.\ Lecture Notes Series
\ \bf 179}.\ Cambridge Univ.\ Press:\ 1992, 106-127.
\vskip .2cm
\item{[CU]} Cukierman,F., Ulmer,D.:\ Curves of genus ten on K3
surfaces.\
To appear in {\it Compositio Mathematica}.
\vskip .2cm
\item{[I1]} Iskovskih,V.A.:\ Fano 3-folds I.\ {\it Math.\ USSR Izv.\
\bf 11},
(1977) 485-527.
\vskip .2cm
\item{[I2]} Iskovskih,V.A.:\ Fano 3-folds II.\ {\it Math.\ USSR Izv.\
\bf 12},
(1978) 479-506.
\vskip .2cm
\item{[KO]} Kobayashi,S., Ochiai,T.:\ Characterizations of complex
 projective
spaces and hyperquadratics.\ {\it J.\ Math.\ Kyoto Univ.\ \bf 13},
(1973)
31-47.
\vskip .2cm
\item{[L]} L'vovskii,S.M.:\ On the extension of varieties defined by
quadratic
equations.\ {\it Math.\ USSR Sbornik \bf 63}, (1989) 305-317.
\vskip .2cm
\item{[M]} Murre,J.P.:\ Classification of Fano threefolds according
to Fano
and Iskovskih.\ In:\ {\it Algebraic Threefolds, Proceedings Varenna
1981.\
Lecture Notes in Math.\ \bf 947}.\ Springer, Berlin-New York:\ 1982,
35-92.
\vskip .2cm
\item{[M1]} Mukai,S.:\ Fano 3-folds.\ In: {\it Complex Projective
Geometry,
Proceedings Trieste-Bergen 1989.\ London Math.\ Soc.\ Lecture Notes
Series \ \bf
179}.\ Cambridge Univ.\ Press:\ 1992, 255-263.
\vskip .2cm
\item{[M2]} Mukai,S.:\ Biregular classification of Fano 3-folds and
Fano
manifolds of coindex 3.\ In: {\it Proc. \ Natl.\ Acad.\ Sci.\ USA
\bf 86}.\
(1989), 3000-3002.
\vskip .2cm
\item{[MM]} Mori,S., Mukai,S.:\ The uniruledness of the moduli space
of curves of
genus 11.\ In:\ {\it Algebraic Geometry, Proceedings Tokyo/Kyoto 1982.\
Lecture
Notes in Math.\ \bf 1016}.\ Springer, Berlin-New York:\ 1983, 334-353.
\vskip .2cm
\item{[W]} Wahl,J.:\ The Jacobian algebra of a graded Gorenstein
singularity.  {\it Duke Math.\ J.\ \bf 55}, (1987) 843-871.
\vskip .2cm
\item{[Z]} Zak,F.L.:\ Some properties of dual varieties and their
application
in projective geometry.\ In:\ {\it Algebraic Geometry, Proceedings
Chicago 1989.
\ Lecture Notes in Math.\ \bf 1479}.\ Springer, Berlin-New York:\
1991, 273-280.

\vskip .5cm
ADDRESSES OF THE AUTHORS:
\vskip .2cm
Ciro Ciliberto, Dipartimento di Matematica, Universit\`a di Roma II,
Tor Vergata,

Viale della Ricerca Scientifica, 00133 Roma, Italy

e-mail: ciliberto@mat.utovrm.it
\vskip .2cm
Angelo Felice Lopez, Dipartimento di Matematica,
Terza Universita' di Roma,

Via Corrado Segre 2, 00146 Roma.

e-mail: lopez@ipvian.bitnet
\vskip .2cm
Rick Miranda, Department of Mathematics, Colorado State University,
Ft. Collins,

CO 80523, USA

e-mail: miranda@riemann.math.colostate.edu

\end